\documentstyle[11pt]{article}
\begin{document}
\author{{Shiqiu Zheng$^{1, 2}$\thanks{Corresponding author, E-mail: shiqiumath@163.com(S. Zheng).}\ , \ \ Shoumei Li$^{1}$\thanks{E-mail: lisma@bjut.edu.cn(S.Li).}}
  \\
\small(1, College of Applied Sciences, Beijing University of Technology, Beijing 100124, China)\\
\small(2, College of Sciences, Hebei United University,  Tangshan 063009, China)\\
}
\date{}
\title{\textbf{Representation for filtration-consistent nonlinear expectations under a general domination condition}\thanks{This work is supported by the National Natural Science Foundation of China (No. 11171010) and the Science and Technology Program of Tangshan (No.
13130203z).}}\maketitle

\textbf{Abstract:}\quad In this paper, we consider filtration-consistent nonlinear expectations which satisfy a general domination condition (dominated by ${\cal{E}}^{\phi}$). We show that this kind of nonlinear expectations can be represented by $g$-expectations defined by the solutions of backward stochastic differential equations, whose generators are independent on $y$ and uniformly continuous in $z$.
\\



\section{Introduction}
The $g$-expectation initiated in Peng [17] in 1997, is a kind of nonlinear expectation defined by the solution of backward stochastic differential equation (BSDE) and can be considered as a nonlinear extension of the Girsanov transformation. The original motivation for studying $g$-expectation comes from the theory of expected utility, which is challenged by the famous Allais paradox and Ellsberg paradox. As a nonlinear expectation, $g$-expectation preserves many properties of the classical linear expectations except the linearity. In particular, it is time-consistent. For properties of $g$-expectation and its applications in utility and risk measures, one can see Briand et al. [1], Chen et al. [2], Chen and Epstein [3], Cohen [4], Coquet et al. [5], Delbaen et al. [6], Jia [10, 11], Jiang [12, 13], Ma and Yao [16], Peng [17, 18, 19], Royer [20] and Rosazza Gianin [21], and among many others.

Time-consistency is one of important properties of $g$-expectation, which will change based on the new observations as time goes on. As a natural extension of $g$-expectation, the notion of filtration-consistent nonlinear expectation is firstly introduced in Coquet et al. [5]. A axiomatic system of this dynamically nonlinear expectation is further introduced in Peng [19]. A very important and interesting result in the Coquet et al. [5] shows that a filtration-consistent nonlinear expectation ${\cal{E}}$ can be represented by a $g$-expectation defined by the solution of a BSDE whose generator $g$ is independent on $y$ and Lipschitz in $z$, when it is translation invariant and satisfies the following domination condition:
$${\cal{E}}[X]-{\cal{E}}[Y]\leq{\cal{E}}^\mu[X-Y],\eqno(1)$$
where ${\cal{E}}^\mu$ is a $g$-expectation defined by the solution of a BSDE whose generator $g=\mu|z|$ for some constant $\mu>0.$

As some extensions of the representation theorem in Coquet et al. [5], Royer [20] obtains a result based on BSDE with jump whose generator $g$ is Lipschitz continuous. Cohen [4] obtains a result based on BSDE in general probability space, whose generator $g$  is also Lipschitz continuous. In fact, the domination conditions in Royer [20] and Cohen [4] are both similar to domination condition (1). Hu et al. [9] obtains a result based on BSDE whose generator $g$  has a quadratic growth, under three domination conditions (see Hu et al. [9, Definition 3.8]) and some other extra conditions. In fact, filtration-consistent nonlinear expectations have a direct correspondence to a fairly large class of risk measures in finance (see Peng [19]). Furthermore, in Hu et al. [9, Page 1519], the authors give the following consequence in finance:
  \begin{center}
  \emph{Any time-consistent risk measure satisfying the required domination condition can be represented by the solution of a simple BSDE!}
  \end{center}
In this topic,  a interesting problem is that can we represent filtration-consistent nonlinear expectation by $g$-expectation under the following domination condition (2)?
$${\cal{E}}[X|{\cal{F}}_t]-{\cal{E}}[Y|{\cal{F}}_t]\leq{\cal{E}}^{\phi}[X-Y|{\cal{F}}_t],\ \ \forall t\in[0,T],\eqno(2)$$
where ${\cal{E}}^\phi$ is a $g$-expectation defined by the solution of a BSDE whose generator $g=\phi(|z|),$ and $\phi(\cdot):{\mathbf{R_+}}\rightarrow{\mathbf{R_+}},$ is subadditive and increasing with $\phi(0)=0$ and has a linear growth. This problem is considered as an unsolved problem by Jia in 2010 (see Jia [11, Remark 4.6]).

In this paper, we answer this problem, using some methods derived from Coquet et al. [5], Hu et al. [9] and Peng [19]. To solve this problem, we will come across some new difficulties, one of which is the most fundamental. That is, the estimation $E|{\cal{E}}^\mu[X-Y]|^2\leq CE|X-Y|^2,$ which holds true for ${\cal{E}}^\mu$, is not always true for ${\cal{E}}^{\phi}$, where $C>0$ is a constant. As a result of this, we can not find a fixed point method can be used directly to solve the BSDE considered in Coquet et al. [5, Theorem 6.1], when such BSDE has a $L^2$  terminal variable and the filtration-consistent expectation ${\cal{E}}$ is dominated by ${\cal{E}}^{\phi}.$ In fact, solving such BSDEs under domination condition (2) is crucial in our paper. Inspired by Hu et al. [9], we use the following strategy. Under domination condition (2), we consider a class of special BSDEs under filtration-consistent expectation ${\cal{E}}$ with bounded terminal variable. Such BSDEs can help us obtain a Doob-Meyer decomposition for ${\cal{E}}$-supermartingale with special construct. Finally, this special Doob-Meyer decomposition is sufficient to establish our representation theorem under domination condition (2). By our representation theorem, we conversely can obtain a existence and uniqueness of BSDEs under ${\cal{E}}$ with $L^2$ terminal variable and a general Doob-Meyer decomposition for ${\cal{E}}$-supermartingale.

This paper is organized as follows. In the next section, we will recall the definitions of $g$-expectation, $g$-martingale and some important results. In section 3, we will recall the definitions of filtration-consistent expectation ${\cal{E}}$, ${\cal{E}}$-martingale and prove some useful properties.  In section 4,  we will give a Doob-Meyer decomposition for ${\cal{E}}$-supermartingale with special construct. In section 5, a representation theorem for filtration-consistent expectation is obtained under domination condition (2).

\section{Preliminaries}
Let $(\Omega,\cal{F},\mathit{P})$ be a complete probability space
carrying a $d$-dimensional standard Brownian motion ${{(B_t)}_{t\geq
0}}$, let $({\cal{F}}_t)_{t\geq 0}$ denote
the natural filtration generated by ${{(B_t)}_{t\geq 0}}$, augmented
by the $\mathit{P}$-null sets of ${\cal{F}}$. Let $|z|$ denote its
Euclidean norm, for $\mathit{z}\in {\mathbf{R}}^d$, $T>0$ be a given time horizon. For stopping times $\tau_1$ and $\tau_2$ satisfying $\tau_1\leq \tau_2,$ Let ${\cal{T}}_{\tau_1,\tau_2}$ be the set of all stopping times $\tau$ satisfying $\tau_1\leq \tau\leq \tau_2.$ For $\tau\in{\cal{T}}_{0,T},$ we define the following usual
spaces:

$L^2({\mathcal {F}}_\tau;{\mathbf{R}}^d)=\{\xi:\ {\cal{F}}_\tau$-measurable
${\mathbf{R}}^d$-valued random variable; ${\mathbf{E}}\left[|\xi|^2\right]<\infty\};$

$L^\infty({\mathcal {F}}_\tau;{\mathbf{R}}^d)=\{\xi:\ {\cal{F}}_\tau$-measurable
${\mathbf{R}}^d$-valued random variable; $\|\xi\|_{L^\infty}=\textrm{esssup}_{\omega\in\Omega}|\xi|<\infty\};$

$L^2_{\cal{F}}(0,\tau;{\mathbf{R}}^d)=\{\psi:\ {\mathbf{R}}^d$-valued predictable
process; $E\left[\int_0^\tau|\psi_t|^2dt\right]
<\infty \};$

$L^\infty_{\cal{F}}(0,\tau;{\mathbf{R}}^d)=\{\psi:\ {\mathbf{R}}^d$-valued predictable
process; $\|\psi\|_{L^\infty_{\cal{F}}}=\textrm{esssup}_{(\omega,t)\in\Omega\times [0,T]}|\psi _t| <\infty\};$

${\mathcal{D}}^2_{\cal{F}}(0,\tau;{\mathbf{R}}^d)=\{\psi:\ $ RCLL
process in $L^2_{\cal{F}}(0,\tau;{\mathbf{R}}^d)$;\ $E
[{\mathrm{sup}}_{0\leq t\leq \tau} |\psi _t|^2] <\infty \}$

${\mathcal{D}}^\infty_{\cal{F}}(0,\tau;{\mathbf{R}}^d)=\{\psi:\ $ RCLL
process in $L^\infty_{\cal{F}}(0,\tau;{\mathbf{R}}^d) \};$

${\mathcal{S}}^2_{\cal{F}}(0,\tau;{\mathbf{R}}^d)=\{\psi:\ $ continuous
process in ${\mathcal{D}}^2_{\cal{F}}(0,\tau;{\mathbf{R}}^d)\};$

${\mathcal{S}}^\infty_{\cal{F}}(0,\tau;{\mathbf{R}}^d)=\{\psi:\ $ continuous
process in ${\mathcal{D}}^\infty_{\cal{F}}(0,\tau;{\mathbf{R}}^d)\}.$\\
Note that when $d=1,$ we always denote $L^2({\mathcal {F}}_\tau;{\mathbf{R}}^d)$ by $L^2({\mathcal {F}}_\tau)$ for convention and use the same treatment for above notations of other spaces.\\

Let us consider a function $g$
$${g}\left( \omega ,t,y,z\right) : \Omega \times [0,T]\times \mathbf{%
R\times R}^{\mathit{d}}\longmapsto \mathbf{R},$$ such that
$\left(g(t,y,z)\right)_{t\in [0,T]}$ is progressively measurable for
each $(y,z)\in\mathbf{%
R\times R}^{\mathit{d}}$. For the function $g$, in this paper, we make the following assumptions:
\begin{itemize}
  \item (A1). There exists a constant
  $K\geq0$ and a continuous function $\phi(\cdot)$, such that $P$-$a.s.,\ \forall t\in[0,T],\
  \forall (y_i,z_i)\in {\mathbf{ R\times R}}^{\mathit{d}},\ \
  (i=1,2):$
  $$|{g}( t,y_{1},z_{1})-{g}( t,y_{2},z_{2}) |\leq K|y_{\mathrm 1}
    -y_{2}|+\phi(|z_{\mathrm1}-z_{2}|),$$
  where $\phi(\cdot):{\mathbf{R_+}}\rightarrow{\mathbf{R_+}},$ is subadditive and increasing with $\phi(0)=0$ and has a linear growth with constant $\nu$, i.e., $\forall x\in {\mathbf{R}}^d, \ \phi(|x|)\leq \nu(|x|+1);$
  \item (A2). $\forall (y,z)\in{\mathbf{ R\times R}}^{\mathit{d}},\ g(t,y,z)\in L^2_{\cal{F}}(0,T);$
  \item (A3). $P$-$a.s.$, $\forall (t,y)\in[0,T]\times{\mathbf{R}},\ g(t,y,0)=0.$
   \item (A1)$^*$. There exists a constant
$\mu\geq0$, such that $P$-$a.s.,\ \forall t\in[0,T],\
\forall (y_i,z_i)\in {\mathbf{ R\times R}}^{\mathit{d}},\ \
(i=1,2):$
$$|{g}( t,y_{1},z_{1})-{g}( t,y_{2},z_{2}) |\leq \mu(|y_{\mathrm 1}
    -y_{2}|+|z_{\mathrm1}-z_{2}|).$$
\end{itemize}

We consider the following BSDEs with parameter $(g,\xi,T):$
$$Y_t=\xi +\int_t^Tg\left(s,Y_s,Z_s\right)
ds-\int_t^TZ_sdB_s,\ \ \forall t\in[0,T].$$
If the generator $g$ satisfies (A1) and (A2), $\xi\in L^2({\mathcal {F}}_T)$, then the BSDE has a unique solution $(Y_t,Z_t)\in{\mathcal{S}}^2_{\cal{F}}(0,T)\times L^2_{\cal{F}}(0,T;{\mathbf{R}}^d)$ (see Jia [10, Theorem 3.3.6], Jia [11, Theorem 2.3] or Fan and Jiang [8, Theorem 2]). Note that since $\phi$ given in (A1) is subadditive and increasing, the BSDE with parameter $(\phi(|z|),\xi,T)$  (resp. $(-\phi(|z|),\xi,T)$) has a unique solution. If $g$ satisfy (A1), (A2) and (A3), a new $g$-expectation and corresponding $g$-martingale are introduced in Jia [10, 11], they are extensions of standard $g$-expectation and $g$-martingale introduced by Peng [17, 18, 19] under (A1)$^*$, (A2) and (A3).\\\\
\
\textbf{Definition 2.1}  Let $g$ satisfy (A1), (A2) and (A3), $\xi\in L^2({\mathcal {F}}_T)$, $(Y_t,Z_t)\in{\mathcal{S}}^2_{\cal{F}}(0,T)\times L^2_{\cal{F}}(0,T;{\mathbf{R}}^d)$ is the solution of BSDE with parameter $(g,\xi,T)$.  The conditional $g$-expectation of $\xi$ is defined by
$${\cal{E}}^g[\xi|{\cal{F}}_t]:=Y_t$$
for $t\in[0,T]$ and $g$-expectation of $\xi$ is defined by
$${\cal{E}}^g[\xi]:=Y_0.$$
\textbf{Definition 2.2} Let $g$ satisfy (A1), (A2) and (A3). A process $Y_t$ with $Y_t\in L^2({\cal{F}}_t)$ for $t\in [0,T],$ is called a $g$-martingale (resp. $g$-supermartingale, $g$-submartingale), if, for each $s\leq t\leq T,$ we have
\begin{center}
${\cal{E}}^g[Y_t|{\cal{F}}_s]=Y_s,$\ \ \ (resp. $\leq,\ \geq$).
\end{center}

Note that we denote ${\cal{E}}^{g}$ by ${\cal{E}}^{\phi}$ (resp. denote ${\cal{E}}^{g}$ by ${\cal{E}}^{-\phi}$), if $g=\phi(|z|)$ (resp. $g=-\phi(|z|)$) for a function $\phi(\cdot)$, and denote ${\cal{E}}^{g}$ by ${\cal{E}}^{\mu}$ (resp. denote ${\cal{E}}^{g}$ by ${\cal{E}}^{-\mu}$),  if $g=\mu|z|$  (resp. $g=-\mu|z|$), for constant $\mu>0$. In fact, following Peng [19], we also can define $g$-martingale (resp. $g$-supermartingale, $g$-submartingale) without (A3), only under (A1) (or (A1)$^*$) and (A2).\\

The following Lemma 2.1 coming from Jia [10, Theorem 3.6.11] is the Doob-Meyer decomposition of $g$-supermartingale under (A1) and (A2). \\\\
\
\textbf{Lemma 2.1} Let $g$ satisfies (A1) and (A2), $Y_t$ is a $g$-supermartingale and has right-continuous path. Then there exists a RCLL process $A_t$, which is increasing with $A_0=0$ and $A_T\in L^2({\cal{F}}_T),$ such that $(Y_t, Z_t)$ is the solution of the following BSDE:
$$Y_t=\xi +\int_t^Tg\left(s,Y_s,Z_s\right)
ds+A_T-A_t-\int_t^TZ_sdB_s,\ \  t\in[0,T].$$

\section{Filtration-consistent nonlinear expectation}
In this section, we will recall the definitions of filtration-consistent expectation ${\cal{E}}$, ${\cal{E}}$-martingale introduced in Peng [19] and prove some important properties which are useful in the proof of our main result. \\\\
\
\textbf{Definition 3.1} Define a system of operators:
$${\cal{E}}[\cdot|{\cal{F}}_t]:\ L^2({\cal{F}}_T)\longrightarrow L^2({\cal{F}}_t), \ t\in[0,T].$$
The operator ${\cal{E}}[\cdot|{\cal{F}}_t]$ is called filtration-consistent condition expectation (${\cal{F}}$-expectation for short), if it satisfies the following aximos:

(i) Monotonicity: ${\cal{E}}[\xi|{\cal{F}}_t]\geq{\cal{E}}[\eta|{\cal{F}}_t], P-a.s., $ if $\xi\geq \eta,\ P-a.s.;$

(ii) Constant preservation: ${\cal{E}}[\xi|{\cal{F}}_t]=\xi, P-a.s., $ if $\xi\in L^2({\mathcal {F}}_t);$

(iii) Consistency: ${\cal{E}}[{\cal{E}}[\xi|{\cal{F}}_t|{\cal{F}}_s]={\cal{E}}[\xi|{\cal{F}}_s], P-a.s.,$ if $s\leq t\leq T;$

(iv) "0-1 Law": ${\cal{E}}[1_A\xi|{\cal{F}}_t]=1_A{\cal{E}}[\xi|{\cal{F}}_t], P-a.s., $ if $A\in{\mathcal {F}}_t.$\\\\
\
\textbf{Definition 3.2} A process $Y_t$ with $Y_t\in L^2({\cal{F}}_t)$ for $t\in [0,T],$ is called a ${\cal{E}}$-martingale (resp. ${\cal{E}}$-supermartingale, ${\cal{E}}$-submartingale), if, for each $s\leq t\leq T,$ we have
\begin{center}
${\cal{E}}[Y_t|{\cal{F}}_s]=Y_s,$\ \ \ (resp. $\leq,\ \geq$).
\end{center}

Note that $g$-expectation defined in Section 2 is an ${\cal{F}}$-expectation (see Jia [11, Theorem 4.3]). Thus the corresponding $g$-martingale (resp. $g$-supermartingale, $g$-submartingale) is also an ${\cal{E}}$-martingale (resp. ${\cal{E}}$-supermartingale, ${\cal{E}}$-submartingale).\\

Now we give some conditions for ${\cal{F}}$-expectation ${\cal{E}}$:
\begin{itemize}
  \item (H1). For each $X,\ Y$ in $L^2({\cal{F}}_T),$ we have
\begin{center}
${\cal{E}}[X|{\cal{F}}_t]-{\cal{E}}[Y|{\cal{F}}_t]\leq{\cal{E}}^{\phi}[X-Y|{\cal{F}}_t], \ \ \forall t \in[0,T],$
\end{center}
where $\phi(\cdot)$ is the function given in (A1).
  \item (H2). (Translation invariance) For each $X$ in $L^2({\cal{F}}_T)$ and $t$ in $[0,T],$ we have,
$${\cal{E}}[X+Y|{\cal{F}}_t]={\cal{E}}[X|{\cal{F}}_t]+Y, \ \ \forall Y\in L^2({\cal{F}}_t)$$
 \item (H1)$^*$. For each $X,\ Y$ in $L^2({\cal{F}}_T),$ we have
\begin{center}
${\cal{E}}[X|{\cal{F}}_t]-{\cal{E}}[Y|{\cal{F}}_t]\leq{\cal{E}}^\mu[X-Y|{\cal{F}}_t], \ \ \forall t \in[0,T],$
\end{center}
where $\mu>0$ is a constant.
\end{itemize}

We list some properties of ${\cal{F}}$-expectation ${\cal{E}},$ which are useful in this paper.\\\\
\
\textbf{Lemma 3.1} Let ${\cal{F}}$-expectation ${\cal{E}}$ satisfy (H1). Then for each $X,\ Y, \ X_n$ in $L^2({\cal{F}}_T),\ n\geq1,$ we have, $\forall t\in[0,T]$,

(i) $-{\cal{E}}^{-\phi}[X|{\cal{F}}_t]={\cal{E}}^{\phi}[-X|{\cal{F}}_t];$

(ii) ${\cal{E}}^{-\phi}[X-Y|{\cal{F}}_t]\leq{\cal{E}}[X|{\cal{F}}_t]-{\cal{E}}[Y|{\cal{F}}_t]\leq{\cal{E}}^{\phi}[X-Y|{\cal{F}}_t];$

(iii) ${\cal{E}}^{-\phi}[X|{\cal{F}}_t]\leq{\cal{E}}[X|{\cal{F}}_t]\leq{\cal{E}}^{\phi}[X|{\cal{F}}_t];$

(iv) $|{\cal{E}}[X|{\cal{F}}_t]-{\cal{E}}[Y|{\cal{F}}_t]|\leq{\cal{E}}^{\phi}[|X-Y||{\cal{F}}_t];$

(v) $\lim_{n\rightarrow\infty}E[|{\cal{E}}[X_n|{\cal{F}}_t]-{\cal{E}}[X|{\cal{F}}_t]|^2]=0,$  if   $\lim_{n\rightarrow\infty}E[|X_n-X|^2]=0;$
\\\\
\
\emph{Proof.} (i) can be checked immediately. By (i) and (H1), we have
 $${\cal{E}}^{-\phi}[X-Y|{\cal{F}}_t]=-{\cal{E}}^{\phi}[Y-X|{\cal{F}}_t]\leq{\cal{E}}[X|{\cal{F}}_t]-{\cal{E}}[Y|{\cal{F}}_t]
 \leq{\cal{E}}^{\phi}[X-Y|{\cal{F}}_t],$$
 then (ii) holds true. (iii) is a consequence of (ii) and "Constant preservation" of ${\cal{E}}$. By (i), (ii) and "Monotonicity" of ${\cal{E}}$, we have
$$-{\cal{E}}^{\phi}[|X-Y||{\cal{F}}_t]={\cal{E}}^{-\phi}[-|X-Y||{\cal{F}}_t]\leq{\cal{E}}[X|{\cal{F}}_t]-{\cal{E}}[Y|{\cal{F}}_t]\leq{\cal{E}}^{\phi}[|X-Y||{\cal{F}}_t]$$
then (iv) holds true. If $\lim_{n\rightarrow\infty}E[|X_n-X|^2]=0,$ by the "Constant preservation" of ${\cal{E}}^{\phi}$ and Jia [11, Theorem 3.11], we can get
$$\lim_{n\rightarrow\infty}E\left[{\cal{E}}^{\phi}[|X_n-X||{\cal{F}}_t]\right]^2=
\lim_{n\rightarrow\infty}E\left[{\cal{E}}^{\phi}[|X_n-X||{\cal{F}}_t]-{\cal{E}}^{\phi}[0|{\cal{F}}_t]\right]^2=0.$$
Then combining above equality and (iv), we obtain (v). \ \ $\Box$ \\\\
\textbf{Remark 3.1} \begin{itemize}
                      \item Let ${\cal{F}}$-expectation ${\cal{E}}$ satisfy (H1). By (H1) and (ii) in Lemma 3.1, for each $X$ in $L^2({\cal{F}}_T)$ and $t$ in $[0,T],$ we have,
$$Y={\cal{E}}^{-\phi}[Y|{\cal{F}}_t]\leq{\cal{E}}[X+Y|{\cal{F}}_t]-{\cal{E}}[X|{\cal{F}}_t]\leq{\cal{E}}^{\phi}[Y|{\cal{F}}_t]=Y,\ \ \forall Y\in L^2({\cal{F}}_t).$$
Then we can get that (H1) implies (H2). Consequently, (H1)$^*$  implies (H2).
                      \item Coquet et al. [5] shows that an ${\cal{F}}$-expectation is a $g$-expectation defined by the solution of a BSDE whose generator $g$ is independent on $y$ and satisfies  (A1)$^*$, (A2) and (A3), when ${\cal{E}}$ satisfies (H2), domination condition (1) and a strict monotonicity condition.  In fact,  By Coquet et al. [5, Lemma 4.3 and Lemma 4.4] and the fact (H1)$^*$  implies (H2), we can easily check that the strict monotonicity condition for ${\cal{E}}$ in Coquet et al. [5] guarantees that (H1)$^*$ is equivalent to (H2) plus domination condition (1).
                    \end{itemize}
\textbf{Lemma 3.2} Let ${\cal{F}}$-expectation ${\cal{E}}$ satisfy (H1). Then for each $X$ in $L^2({\cal{F}}_T),$ ${\cal{E}}[X|{\cal{F}}_t]$ admits a RCLL version.\\\\
\
\emph{Proof.} Since $\phi(\cdot)$ has a linear growth, by Lepeltier and San Martin [15, Lemme 1], we can find a function $\varphi(z): {\mathbf{R}}_+\mapsto \mathbf{R},$ which satisfies (A1)$^*$ and (A2), such that for each $z\in{\mathbf{R}}^{\mathit{d}}$, $\varphi(|z|)\leq {-\phi(|z|)}.$ By (iii) in Lemma 3.1 and comparison theorem (see Jia [11, Theorem 3.1]), we have $${\cal{E}}[X|{\cal{F}}_t]\geq {\cal{E}}^{-\phi}[X|{\cal{F}}_t]\geq {\cal{E}}^{\varphi}[X|{\cal{F}}_t].$$
Consequently, we can easily check that ${\cal{E}}[X|{\cal{F}}_t]$ is a $\varphi$-supermartingale. By Peng [19, Theorem 3.7], we get that for a denumerable dense subset ${{D}}$ of $[0,T]$, almost all $\omega$ and all $t\in[0,T],$ we have $\lim_{s\in{{D}},\ s\searrow t}{\cal{E}}[X|{\cal{F}}_s]$ and $\lim_{s\in{{D}},\ s\nearrow t}{\cal{E}}[X|{\cal{F}}_s]$ both exist and are finite. For each $t\in[0,T]$, we set
$$Y_t:=\lim_{s\in{{D}},\ s\searrow t}{\cal{E}}[X|{\cal{F}}_s],\eqno(3)$$
then $Y_t$ is RCLL.  By (iv) in Lemma 3.1 and "Constant preservation" of ${\cal{E}}$, we have
$$|{\cal{E}}[X|{\cal{F}}_t]|\leq{\cal{E}}^{\phi}[|X||{\cal{F}}_t].\eqno(4)$$
By Jia [11, Theorem 2.3], we also have
$$ E[\sup_{t\in[0,T]}|{\cal{E}}^{\phi}[|X||{\cal{F}}_t]|^2]<+\infty.\eqno(5)$$
By (3)-(5) and Lebesgue dominated convergence theorem, we have
$$\lim_{s\in{{D}},\ s\searrow t}{\cal{E}}[X|{\cal{F}}_s]=Y_t, \ \ \forall t\in[0,T]. \eqno(6)$$
in $L^2({\cal{F}}_T)$ sense.
By the "Constant preservation" of ${\cal{E}}$, we have
$$
  {\cal{E}}[{\cal{E}}[X|{\cal{F}}_s]|{\cal{F}}_t]-Y_t= {\cal{E}}[{\cal{E}}[X|{\cal{F}}_s]|{\cal{F}}_t]-{\cal{E}}[Y_t|{\cal{F}}_t].
$$
Then by (6) and (v) in Lemma 3.1,  we get
$$\lim_{s\in{{D}},\ s\searrow t}{\cal{E}}[{\cal{E}}[X|{\cal{F}}_s]|{\cal{F}}_t]=Y_t,\ \ \forall t\in[0,T], \eqno(7)$$
in $L^2({\cal{F}}_T)$ sense. On the other hand, by "Consistency" of ${\cal{E}}$, we have
$$\lim_{s\in{{D}},\ s\searrow t}{\cal{E}}[{\cal{E}}[X|{\cal{F}}_s]|{\cal{F}}_t]={\cal{E}}[X|{\cal{F}}_t],\ \ \forall t\in[0,T].\eqno(8)$$
By (7) and (8), we have $P-a.s.,\ {\cal{E}}[X|{\cal{F}}_t]=Y_t.$ The proof is complete.\ \  $\Box$ \\

Note that, in the sequel, we always take the RCLL version of ${\cal{E}}[X|{\cal{F}}_t]$. If for each $t\in[0,T]$, $Y_t={\cal{E}}[X|{\cal{F}}_t]$, then for each stopping time $\sigma\in{\cal{T}}_{0,T},$ we set
${\cal{E}}[X|{\cal{F}}_\sigma]:=Y_\sigma.$ Then we have the following optional stopping theorem, which can be proved by Lemma 3.1 and the same arguments as Peng [19, Theorem 7.4], directly. We omit its proof here.\\\\
\textbf{Lemma 3.3} Let ${\cal{F}}$-expectation ${\cal{E}}$ satisfy (H1). If $Y_t\in {\cal{D}}^2_{{\cal{F}}}(0,T)$ is an ${\cal{E}}$-supermartingale (resp. ${\cal{E}}$-submartingale), then for each stopping time $\sigma,\ \tau\in{\cal{T}}_{0,T},$ we have
$${\cal{E}}[Y_\tau|{\cal{F}}_\sigma]\leq Y_{\sigma\wedge\tau},\ (\textrm{resp.} \geq Y_{\sigma\wedge\tau}),\ \ P-a.s.$$

The following Lemma 3.4 can be considered as a representation theorem for ${\cal{E}}$-martingale.\\\\
\textbf{Lemma 3.4} Let ${\cal{F}}$-expectation ${\cal{E}}$ satisfy (H1) and for each $X\in L^2({\cal{F}}_T),$ set $$y_t^X:={\cal{E}}[X|{\cal{F}}_t],\ \forall t\in[0,T].$$ Then there exists a pair $(g_t^X, Z_t^X)$ in $L^2_{{\cal{F}}}(0,T)\times L^2_{{\cal{F}}}(0,T;{\textbf{R}}^d)$ such that $$|g_t^X|\leq \phi(|Z_t^X|),\eqno(9)$$ and
$$y_t^X=X+\int_t^Tg_s^Xds-\int_t^TZ_s^XdB_s.\eqno(10)$$
Moreover, for $Y$ in $L^2({\cal{F}}_T),$ we have $$|g_t^X-g_t^Y|\leq  \phi(|Z_t^X-Z_t^Y|).\eqno(11)$$
\emph{Proof.} By (iii) in Lemma 3.1 and Lemma 3.2, we can easily check that ${\cal{E}}[X|{\cal{F}}_t], t\in[0,T]$ is a right-continuous ${-\phi}$-supermartingale (resp. ${\phi}$-submartingale). Then by Lemma 2.1, we get that there exists $(Z^{X,-\phi}_t,A^{X,-\phi}_t)$ (resp. $(Z^{X,\phi}_t,A^{X,\phi}_t)$) in $L^2_{{\cal{F}}}(0,T;{\textbf{R}}^d)\times D^2_{{\cal{F}}}(0,T)$ with $A^{X,-\phi}_t$ (resp. $A^{X,\phi}_t$) RCLL, increasing and $A^{X,-\phi}_0=0$ (resp. $A^{X,\phi}_0=0$), such that
$$y_t^X=X-\int_t^T\phi(|Z_s^{X,-\phi}|)ds+A^{X,-\phi}_T-A^{X,-\phi}_t-\int_t^TZ_s^{X,-\phi}dB_s,\eqno(12)$$
and
$$y_t^X=X+\int_t^T\phi(|Z_s^{X,\phi}|)ds-A^{X,\phi}_T+A^{X,\phi}_t-\int_t^TZ_s^{X,\phi}dB_s.\eqno(13)$$
Comparing the martingale parts and the bounded variation parts of (12) and (13), we get
\begin{eqnarray*}
  Z_s^{X,-\phi}&\equiv& Z_s^{X,\phi},\\
-\phi(|Z_s^{X,-\phi}|)ds+dA^{X,-\phi}_s&\equiv&\phi(|Z_s^{X,\phi}|)ds-dA^{X,\phi}_s.
\end{eqnarray*}
Thus we have
$$2\phi(|Z_s^{X,\phi}|)ds\equiv dA^{X,-\phi}_s+dA^{X,\phi}_s.$$
Consequently, $A^{X,\phi}_s$ and $A^{X,-\phi}_s$ are both absolutely continuous. Thus there exist $a^{X,\phi}_s\geq0$ and $a^{X,-\phi}_s\geq0$ such that
$$dA^{X,\phi}_s=a^{X,\phi}_sds,\ \ dA^{X,-\phi}_s=a^{X,-\phi}_sds.$$
Then we have
$$a^{X,\phi}_s+a^{X,-\phi}_s\equiv 2\phi(|Z_s^{X,\phi}|).$$
By setting
$$Z^X:=Z_s^{X,\phi},\ \ \ g^X:=\phi(|Z_s^{X}|)-a^{X}_sds,$$
we get (9) and (10).
By (10), for $Y$ in $L^2({\cal{F}}_T),$ there exists a pair $(g_t^Y, Z_t^Y)$ in $L^2_{{\cal{F}}}(0,T)\times L^2_{{\cal{F}}}(0,T;{\textbf{R}}^d)$ such that
$$y_t^Y=Y+\int_t^Tg_s^Yds-\int_t^TZ_s^YdB_s.\eqno(14)$$
By (10) and (14), we have
$$y_t^X-y_t^Y=X-Y+\int_t^T(g_s^X-g_s^Y)ds-\int_t^T(Z_s^X-Z_s^Y)dB_s.\eqno(15)$$
By (ii) in Lemma 3.1, we have, for each $s\leq t\leq T,$
$${\cal{E}}^{-\phi}[y_t^X-y_t^Y|{\cal{F}}_s]\leq{\cal{E}}[y_t^X|{\cal{F}}_s]-{\cal{E}}[y_t^Y|{\cal{F}}_s]
=y_s^X-y_s^Y={\cal{E}}[y_t^X|{\cal{F}}_s]-{\cal{E}}[y_t^Y|{\cal{F}}_s]\leq {\cal{E}}^{\phi}[y_t^X-y_t^Y|{\cal{F}}_s].$$
Thus $y_t^X-y_t^Y$ is a ${\phi}$-submaringale and a ${-\phi}$-supermaringale. Then by Lemma 2.1 again, there exist $(Z^1_t,A^1_t)$ (resp. $(Z^2_t,A^2_t)$) in $L^2_{{\cal{F}}}(0,T;{\textbf{R}}^d)\times D^2_{{\cal{F}}}(0,T)$ with $A^1_t$ (resp. $A^2_t$) RCLL, increasing and $A^1_0=0$ (resp. $A^2_0=0$), such that
$$y_t^X-y_t^Y=X-Y+\int_t^T\phi(|Z_s^1|)ds-A_T^1+A_t^1-\int_t^TZ_s^1dB_s.\eqno(16)$$
$$y_t^X-y_t^Y=X-Y-\int_t^T\phi(|Z_s^2|)ds+A_T^2-A_t^2-\int_t^TZ_s^2dB_s.\eqno(17)$$
Comparing the martingale parts and the bounded variation parts of (15) and (16), we get
\begin{eqnarray*}
  Z_s^X-Z_s^Y&\equiv& Z_s^1,\\
(g_s^X-g_s^Y)ds&\equiv&\phi(|Z_s^1|)ds-dA_s^1.
\end{eqnarray*}
Then we have
$$g_t^X-g_t^Y\leq \phi(|Z_t^X-Z_t^Y|).\eqno(18)$$
Comparing the martingale parts and the bounded variation parts of (15) and (17), we get
\begin{eqnarray*}
  Z_s^X-Z_s^Y&\equiv& Z_s^2,\\
(g_s^X-g_s^Y)ds&\equiv&-\phi(|Z_s^2|)ds+dA_s^2.
\end{eqnarray*}
Then we have
$$g_t^X-g_t^Y\geq -\phi(|Z_t^X-Z_t^Y|).\eqno(19)$$
Thus by (18) and (19), we obtained (11). The proof is completed. \ \  $\Box$\\\\
\
\textbf{Remark 3.2} Let ${\cal{F}}$-expectation ${\cal{E}}$ satisfy (H1), $X\in L^2({\cal{F}}_T)$ and $\eta_t\in L^2_{{\cal{F}}}(0,T).$ By (10), ${\cal{E}}[X|{\cal{F}}_t]$ is a continuous process. By (H2),  we can further get ${\cal{E}}[X+\int_t^T\eta_sds|{\cal{F}}_t]$ is continuous, from the fact that $${\cal{E}}[X+\int_t^T\eta_sds|{\cal{F}}_t]={\cal{E}}[X+\int_0^T\eta_sds|{\cal{F}}_t]-\int_0^t\eta_sds.$$

The following Lemma 3.5 describes a property of ${\cal{F}}$-expectation, which plays an important role in this paper.\\\\
\textbf{Lemma 3.5} Let ${\cal{F}}$-expectation ${\cal{E}}$ satisfy (H1). Then for each $z\in \textbf{R}^d$ and each $X\in L^\infty({\cal{F}}_T),$ there exists a process $\eta_t\in {\cal{S}}^\infty_{{\cal{F}}}(0,T)$ such that for each $t\in[0,T],$
$${\cal{E}}[X+zB_T|{\cal{F}}_t]=\eta_t+zB_t,\ \ P-a.s.$$
\emph{Proof.} For $z\in \textbf{R}^d$ and $X\in L^\infty({\cal{F}}_T),$ by Lemma 3.4, there exists a pair $(g_t, Z_t)$ in $L^2_{{\cal{F}}}(0,T)\times L^2_{{\cal{F}}}(0,T;{\textbf{R}}^d)$  such that
$${\cal{E}}[X+zB_T|{\cal{F}}_t]=X+zB_T+\int_t^Tg_sds-\int_t^TZ_sdB_s.$$
Then we have
$${\cal{E}}[X+zB_T|{\cal{F}}_t]-zB_t=X+\int_t^Tg_sds-\int_t^T(Z_s-z)dB_s.\eqno(20)$$
Set $$(\tilde{Y}_t,\tilde{Z}_t):=({\cal{E}}[X+zB_T|{\cal{F}}_t]-zB_t,Z_t-z).\eqno(21)$$
Then by (20) and (21), $(\tilde{Y}_t,\tilde{Z}_t)$ is the unique solution of BSDE with parameter $(g_s,X,T).$
Now, we consider the following two BSDEs:
$$\overline{{Y}}_t=X+\int_t^T(\phi (|\overline{{Z}}_s|)+\phi (|z|))ds-\int_t^T\overline{Z}_sdB_s.\eqno(22)$$
$$\underline{{Y}}_t=X-\int_t^T(\phi (|\underline{{Z}}_s|)+\phi (|z|))ds-\int_t^T\underline{Z}_sdB_s.\eqno(23)$$
By (9), (21) and the fact that $\phi$ is increasing and subadditive, we have
$$|g_t|\leq \phi(|Z_t|)=\phi(|\tilde{Z}_t+z|)\leq  \phi(|\tilde{Z}_t|)+\phi(|z|),\eqno(24)$$
Then by (20)-(24) and comparison theorem (see Jia [11, Theorem 3.1]), we can get that for each $t\in[0,T],$
$$\underline{{Y}}_t\leq\tilde{Y}_t\leq\overline{{Y}}_t, \ \ P-a.s.\eqno(25)$$
Since $X\in L^\infty({\cal{F}}_T)$ and $\phi$ has a linear growth, then applying the boundness of solution of quadratic BSDEs (whose generator has a quadratic growth in $z$) with bounded terminal variable (see Kobylanski [14, Theorem
2.3]) to (22) and (23), we get that $\overline{{Y}}_t$ and $\underline{{Y}}_t$ both belong to ${\cal{S}}^\infty_{{\cal{F}}}(0,T).$ Then by (25), we have
$$\|\tilde{Y}_t\|_{L^\infty_{\cal{F}}}\leq\|\underline{{Y}}_t\|_{L^\infty_{\cal{F}}}\vee\|\overline{{Y}}_t\|_{L^\infty_{\cal{F}}}< \infty.$$
From (21) and the above inequality, the proof is complete. \ \  $\Box$\\\\
\
\textbf{Lemma 3.6}  Let ${\cal{F}}$-expectation ${\cal{E}}$ satisfy (H1), $\sigma\in{\cal{T}}_{0,T}$ and $X,\ Y\in L^2({\cal{F}}_T).$ Then we have
$$1_A{\cal{E}}[X+Y|{\cal{F}}_\sigma]=1_A{\cal{E}}[1_AX+Y|{\cal{F}}_\sigma],\ \ P-a.s. \ \ \forall A\in {\cal{F}}_{\sigma}.$$
\emph{Proof.} The proof can be completed by Lemma 3.2 and the similar argument as Hu et al. [9, Proposition 4.2(iii)]. We omit it here.\ \  $\Box$\\\\
\
\textbf{Lemma 3.7}  Let ${\cal{F}}$-expectation ${\cal{E}}$ satisfy (H1), $\sigma\in{\cal{T}}_{0,T}$ and $X\in L^2({\cal{F}}_T).$ Then we have
$${\cal{E}}[X+Y|{\cal{F}}_\sigma]={\cal{E}}[X|{\cal{F}}_\sigma]+Y,\ \ P-a.s.,\ \ \forall Y\in L^2({\cal{F}}_\sigma).$$
\emph{Proof.} The proof can be completed by (H2), Lemma 3.2 and the similar argument as Hu et al. [9, Proposition 4.2(iv)]. We omit it here.\ \  $\Box$

\section{Doob-Meyer decomposition of ${\cal{E}}$-supermartingale}
In this section, we will study the Doob-Meyer decomposition of ${\cal{E}}$-supermartingale. Firstly, we consider a BSDE under ${\cal{F}}$-expectation ${\cal{E}}.$\\

Given a function $f: \Omega \times [0,T]\times \mathbf{R}\longmapsto \mathbf{R},$ in this paper, we always suppose $f$ satisfy the following Lipschitz condition.
$$\exists \lambda\geq0,\ s.t.\ |f(t,y_1)-f(t,y_2)|\leq \lambda|y_1-y_2|, \ \forall y_1,\ y_2\in{\mathbf{R}},\ \forall t\in[0,T].$$

Now, we consider the following BSDE under ${\cal{F}}$-expectation ${\cal{E}}:$
$$y_t+zB_t={\cal{E}}\left[X+zB_T+\int_t^Tf(s,y_s)ds|{\cal{F}}_t\right],\ \ t\in[0,T]$$
which has been studied in Coquet et al. [5] for the case $z=0$, $X\in L^2({\cal{F}}_T)$ and $f(\cdot,0)\in L^2_{{\cal{F}}}(0,T),$ and in Hu et al. [9] for the case $z\in\textbf{R}^d$, $X\in L^\infty({\cal{F}}_T)$ and $f(\cdot,0)\in L^\infty_{{\cal{F}}}(0,T).$ We denote this BSDE by ${\cal{E}}(f,T,X,z).$ The following Theorem 4.1 shows that it has a unique solution under (H1). \\\\
\textbf{Theorem 4.1} Let ${\cal{F}}$-expectation ${\cal{E}}$ satisfy (H1), $z\in\textbf{R}^d,$ $X\in L^\infty({\cal{F}}_T)$ and $f(\cdot,0)\in L^\infty_{{\cal{F}}}(0,T).$ Then ${\cal{E}}(f,T,X,z)$ has a unique solution $y_t\in {\cal{S}}_{\cal{F}}^\infty (0,T).$ \\\\
\emph{Proof.} For $y_t\in {\cal{S}}^\infty_{{\cal{F}}}(0,T),$ set $$I(y_t):={\cal{E}}\left[X+zB_T+\int_t^Tf(s,y_s)ds|{\cal{F}}_t\right]-zB_t,$$
By (H2), we have
$$I(y_t)={\cal{E}}\left[X+zB_T+\int_0^Tf(s,y_s)ds|{\cal{F}}_t\right]-zB_t-\int_0^tf(s,y_s)ds,\eqno(26)$$
Since $f$ satisfies Lipschitz condition, $y_t\in {\cal{S}}^\infty_{{\cal{F}}}(0,T)$ and $f(\cdot,0)\in L^\infty_{{\cal{F}}}(0,T)$, thus, we have
\begin{eqnarray*}
\left\|\int_0^tf(s,y_s)ds\right\|_{L^\infty_{\cal{F}}}\leq\int_0^T\|f(s,y_s)\|_{L^\infty_{\cal{F}}} ds&\leq&\int_0^T\|f(s,0)\|_{L^\infty_{\cal{F}}}ds+\lambda\int_0^T\|y_s\|_{L^\infty_{\cal{F}}}ds\\
&=&T\|f(s,0)\|_{L^\infty_{\cal{F}}}+\lambda T\|y_s\|_{L^\infty_{\cal{F}}}\\
&<&\infty.
\end{eqnarray*}
With the help of $X\in L^\infty({\cal{F}}_T),$ the above inequality, Remark 3.2 and Lemma 3.5, we can get $I(y_t)\in {\cal{S}}_{\cal{F}}^\infty (0,T) $ from (26), Thus
$$I(\cdot):{\cal{S}}_{\cal{F}}^\infty (0,T)\longmapsto {\cal{S}}_{\cal{F}}^\infty (0,T).$$
By (iv) in Lemma 3.1, "Monotonicity" and "Constant preservation" of ${\cal{E}}^{\phi},$ for each $t\in[0,T]$ and $y_t^1,\ y_t^2\in {\cal{S}}_{\cal{F}}^\infty (0,T),$ we have
\begin{eqnarray*}
  |I(y_t^1)-I(y_t^2)|&=&\left|{\cal{E}}\left[X+zB_T+\int_t^Tf(s,y_s^1)ds|{\cal{F}}_t\right]
  -{\cal{E}}\left[X+zB_T+\int_t^Tf(s,y_s^2)ds|{\cal{F}}_t\right]\right|\\
  &\leq&{\cal{E}}^{\phi}\left[\left|\int_t^Tf(s,y_s^1)ds
  -\int_t^Tf(s,y_s^2)ds\right||{\cal{F}}_t\right]\\
  &\leq&{\cal{E}}^{\phi}\left[\int_t^T\left|f(s,y_s^1)-f(s,y_s^2)\right|ds|{\cal{F}}_t\right]\\
  &\leq&\lambda T\|y_t^1-y_t^2\|_{L^\infty_{\cal{F}}}
\end{eqnarray*}

\textbf{Case 1: } $T\leq\frac{1}{2\lambda}.$

In this case, we have $\|I(y_t^1)-I(y_t^2)\|_{L^\infty_{\cal{F}}}\leq\frac{1}{2}\|y_t^1-y_t^1\|_{L^\infty_{\cal{F}}}.$ Thus $I(\cdot)$ is a strict contraction. The proof is complete.

\textbf{Case 2: } $T>\frac{1}{2\lambda}.$

In this case, we can complete the proof using a "patching-up" method given in Hu et al. [9, Proposition 4.4]. We omit it here.\ \  $\Box$\\\\
\
\textbf{Remark 4.1 } Let ${\cal{F}}$-expectation ${\cal{E}}$ satisfy (H1) and $y_t$ is the solution of ${\cal{E}}(f,T,X,z),$ then by (H2), we can get the process $y_t+zB_t+\int_0^tf(s,y_s)ds$ is an ${\cal{E}}$-martingale.\\\\
\
\textbf{Theorem 4.2} Let ${\cal{F}}$-expectation ${\cal{E}}$ satisfy (H1), $z\in \textbf{R}^d$, $X\in L^\infty({\cal{F}}_T)$, $f(\cdot,0)\in L^\infty_{{\cal{F}}}(0,T),$ $y_t$ is the solution of ${\cal{E}}(f,T,X,z)$ and $\bar{y}_t$ is the solution of the following ${\cal{E}}(f+\eta_t,T,\bar{X},z)$:
$$\bar{y}_t+zB_t={\cal{E}}\left[\bar{X}+zB_T+\int_t^T(f(s,\bar{y}_s)+\eta_s)ds|{\cal{F}}_t\right],\ \  t\in[0,T],\eqno(27)$$
where $\bar{X}$ in $L^\infty({\cal{F}}_T)$ and $\eta_t\in L_{\cal{F}}^\infty (0,T)$ satisfy
$$\bar{X}\geq X,\ \ \ \eta_t\geq0,\ \ dP\times dt-a.e.$$
Then $\forall t\in[0,T],$ we have
$$\bar{y}_t\geq y_t,\ \ P-a.s.\eqno(28)$$
\emph{Proof.} \textbf{Case 1:} $\eta_t\equiv0.$

For constant $n\geq1$, we define the stopping time
$$\tau_1^n:=\inf\left\{t\geq0;\ \bar{y}_t\leq y_t-\frac{1}{n}\right\}\wedge T.$$
Clearly, if (28) is not true, then there exists a integer $k\geq1$ such that $P(\{\tau_1^k<T\})>0.$
By setting $A:=\{\tau_1^k<T\}$ and the continuity of $\bar{y}_t$ and $y_t$, and $\bar{y}_T=\bar{X}\geq X={y}_T,$ we have
$$P(A)>0\ \ \textrm{and} \ \ A=\left\{\bar{y}_{\tau_1^k}\leq y_{\tau_1^k}-\frac{1}{k}\right\}.\eqno(29)$$
Now, we define the stopping time
$$\tau_2:=\inf\{t\geq\tau_1^k;\ \bar{y}_t\geq y_t\}\wedge T.$$
By the continuity of $\bar{y}_t$ and $y_t$, and $\bar{y}_T=\bar{X}\geq X={y}_T,$ we have
$$1_A\bar{y}_{\tau_2}=1_A{y}_{\tau_2}.$$
Clearly, $A\in{\cal{F}}_{\tau_1^k}.$ Then, for each stopping time $\tau\in{\cal{T}}_{\tau_1^k,\tau_2},$ we have
\begin{eqnarray*}
  &&1_A{\cal{E}}\left[1_A\bar{y}_{\tau_2}+zB_{\tau_2}+\int_{\tau}^{{\tau_2}}1_Af(s,1_A\bar{y}_s)ds|{\cal{F}}_{\tau}\right]\\
  &=& 1_A{\cal{E}}\left[1_A\bar{y}_{\tau_2}+zB_{\tau_2}+1_A\int_{\tau}^{{\tau_2}}f(s,\bar{y}_s)ds|{\cal{F}}_{\tau}\right]\\
  &=& 1_A{\cal{E}}\left[\bar{y}_{\tau_2}+zB_{\tau_2}+\int_{\tau}^{{\tau_2}}f(s,\bar{y}_s)ds|{\cal{F}}_{\tau}\right]\\
  &=& 1_A{\cal{E}}\left[\bar{y}_{\tau_2}+zB_{\tau_2}+\int_{0}^{{\tau_2}}f(s,\bar{y}_s)ds|{\cal{F}}_{\tau}\right]
   -1_A\int_{0}^{{\tau}}f(s,\bar{y}_s)ds\\
  &=& 1_A\left(\bar{y}_{\tau}+zB_{\tau}+\int_{0}^{{\tau}}f(s,\bar{y}_s)ds\right)
   -1_A\int_{0}^{{\tau}}f(s,\bar{y}_s)ds\\
  &=& 1_A(\bar{y}_{\tau}+zB_{\tau})
\end{eqnarray*}
In the above, the second equality is due to Lemma 3.6, the third equality is due to Lemma 3.7, the fourth equality is due to Remark 4.1 and Lemma 3.3.

By the same argument as above, for each stopping time $\tau\in{\cal{T}}_{\tau_1^k,\tau_2},$ we have
$$1_A{\cal{E}}\left[1_Ay_{\tau_2}+zB_{\tau_2}+\int_{\tau}^{\tau_2}1_Af(s,1_Ay_s)ds
|{\cal{F}}_{\tau}\right]= 1_A(y_{\tau}+zB_{\tau}).$$
For each $t\in[0,T]$, we set stopping time $\hat{t}:=(t\vee\tau_1^k)\wedge\tau_2.$ Then by above three equalities, (iv) in Lemma 3.1, "Monotonicity" and "Constant preservation" of ${\cal{E}}^{\phi},$ we can get,
\begin{eqnarray*}
  |1_A\bar{y}_{\hat{t}}-1_Ay_{\hat{t}}|&=&|1_A(\bar{y}_{\hat{t}}+zB_{\hat{t}})-1_A(y_{\hat{t}}+zB_{\hat{t}})|\\
  &\leq&1_A{\cal{E}}^{\phi}\left[\left|\int_{\hat{t}}^{{\tau_2}}1_Af(s,1_A\bar{y}_s)ds
  -\int_{\hat{t}}^{\tau_2}1_Af(s,1_Ay_s)ds\right||{\cal{F}}_{\hat{t}}\right]\\
  &\leq&{\cal{E}}^{\phi}\left[\int_{\hat{t}}^{{\tau_2}}\left|f(s,1_A\bar{y}_s)-f(s,1_Ay_s)\right|ds|{\cal{F}}_{\hat{t}}\right]\\
  &\leq&{\cal{E}}^{\phi}\left[\left\|\int_{\hat{t}}^T\lambda|1_A\bar{y}_s-1_Ay_s|ds\right\|_{L^\infty}|{\cal{F}}_{\hat{t}}\right]\\
  &\leq&\left\|\int_t^T\lambda|1_A\bar{y}_{\hat{s}}-1_Ay_{\hat{s}}|ds\right\|_{L^\infty}\\
  &\leq&\lambda\int_t^T\left\|1_A\bar{y}_{\hat{s}}-1_Ay_{\hat{s}}\right\|_{L^\infty} ds.
\end{eqnarray*}
Then we have, for each $t\in[0,T]$,
$$\|1_A\bar{y}_{\hat{t}}-1_Ay_{\hat{t}}\|_{L^\infty} \leq\lambda\int_t^T\left\|1_A\bar{y}_{\hat{s}}-1_Ay_{\hat{s}}\right\|_{L^\infty} ds.$$
By Gronwall inequality, we have for each $t\in[0,T]$, $1_A\bar{y}_{\hat{t}}=1_Ay_{\hat{t}},\ P-a.s.$ By setting $t=0,$ we have
$$1_A\bar{y}_{\tau_1^k}=1_Ay_{\tau_1^k},\ \ P-a.s.$$
which contradicts (29). Thus (28) holds true.

\textbf{Case 2:} $\eta_t\equiv 0$ is not true.

For $n\geq1,$  set $t^n_i:=\frac{i}{n}T,\ 1\leq i\leq n.$ As in Coquet et al. [5] and Hu et al. [9], we define the following BSDEs recursively
$$y_t^{i,n}+zB_t={\cal{E}}\left[\left(X^n_i+\int_{t^n_{i-1}}^{t^n_i}\eta_sds\right)+zB_{t^n_i}+\int_{t}^{t^n_i}f(s,y_s^{i,n})ds|{\cal{F}}_t\right],\ t\in[0,t^n_i],$$
where $X^n_n=\bar{X}$ and $X^n_{i-1}=y_{t^n_{i-1}}^{i,n},$ for $1\leq i\leq n.$ By the result of Case 1, we have $y_t^{n,n}\geq y_t,\ t\in[t^n_{n-1},{t^n_n}].$ Thus $X^n_{n-1}+\int_{t^n_{n-2}}^{t^n_{n-1}}\eta_sds\geq y_{t^n_{n-1}}.$ then by the result of Case 1 again, we have $y_t^{n-1,n}\geq y_t,\ t\in[{t^n_{n-2}},{t^n_{n-1}}].$ Similarly, we also have $y_t^{i,n}\geq y_t,\ t\in[{t^n_{i-1}},{t^n_{i}}],\ 1\leq i\leq n-2.$ We define the process $y_t^n=y_t^{i,n},\ t\in[{t^n_{i-1}},{t^n_i}),\ 1\leq i\leq n,\ y_T^n=\bar{X}.$ Then we can check that
$$y_t^{n}+zB_t={\cal{E}}\left[\bar{X}+\int_{t^n_{i-1}}^{T}\eta_sds+zB_{T}+\int_{t}^{T}f(s,y_s^n)ds|{\cal{F}}_t\right],\ t\in[{t^n_{i-1}},{t^n_i}),\ \ 1\leq i\leq n.\eqno(30)$$
By (27), (30), (iv) in Lemma 3.1, "Monotonicity" and "Constant preservation" of ${\cal{E}}^{\phi},$ for $t\in[{t^n_{i-1}},{t^n_i}),\ 1\leq i\leq n,$ we have
\begin{eqnarray*}
  &&|y_t^{n}-\bar{y}_t|\\&=&\left|{\cal{E}}\left[\bar{X}+\int_{t^n_{i-1}}^{T}\eta_sds
  +zB_T+\int_{t}^{T}f(s,y_s^n)ds|{\cal{F}}_t\right]-{\cal{E}}\left[\bar{X}+zB_T
  +\int_t^T(f(s,\bar{y}_s)+\eta_s)ds|{\cal{F}}_t\right]\right|\\
  &\leq&{\cal{E}}^{\phi}\left[\left|\int_{t^n_{i-1}}^{t}\eta_sds
  +\int_{t}^{T}(f(s,y_s^n)-f(s,\bar{y}_s))ds\right||{\cal{F}}_t\right]\\
   &\leq&{\cal{E}}^{\phi}\left[\int_{t^n_{i-1}}^{t}|\eta_s|ds
  +\int_{t}^{T}|f(s,y_s^n)-f(s,\bar{y}_s)|ds|{\cal{F}}_t\right]\\
 &\leq&\int_{t^n_{i-1}}^{t}\left\|\eta_s\right\|_{L^\infty} ds
  +\int_{t}^{T}\left\|(f(s,y_s^n)-f(s,\bar{y}_s))\right\|_{L^\infty} ds\\
  &\leq&\frac{T}{n}\left\|\eta_s\right\|_{L_{\cal{F}}^\infty}+\lambda\int_{t}^T\left\|\bar{y}_s-y_s^n\right\|_{L^\infty} ds.
\end{eqnarray*}
By Gronwall inequality, we get for $t\in[0,T],$ $y_t^{n}\rightarrow\bar{y}_t$ in $L^\infty({\cal{F}}_t)$ sense, as $n\rightarrow\infty.$ Consequently, $\forall t\in[0,T],\ \bar{y}_t\geq y_t,\  P-a.s.$ The proof is complete.\ \  $\Box$\\

Now, we give the following Doob-Meyer decomposition of ${\cal{E}}$-supermartingale.\\\\
\textbf{Theorem 4.3}  Let ${\cal{F}}$-expectation ${\cal{E}}$ satisfy (H1), $z\in \textbf{R}^d$, $Y_t\in S^\infty_{{\cal{F}}}(0,T)$  and $Y_t+zB_t$ is an ${\cal{E}}$-supermartingale, then there exists a process $A_t\in S^2_{{\cal{F}}}(0,T)$, which is increasing with $A_0=0$ such that $\forall t\in [0,T],$
\begin{center}
${\cal{E}}[Y_T+zB_T+A_T|{\cal{F}}_t]=Y_t+zB_t+A_t,\ \ P-a.s.$
\end{center}
\emph{Proof.}  For $n\geq1,$ we consider the following BSDEs under ${\cal{F}}$-expectation ${\cal{E}}$:
$$y_t^n+zB_t={\cal{E}}\left[Y_T+zB_T+\int_t^Tn(Y_s-y_t^n)ds|{\cal{F}}_t\right],\ \ \ t\in[0,T].\eqno(31)$$
By Theorem 4.1, the above BSDE (31) has a unique solution $y_t^n\in S^\infty_{{\cal{F}}}(0,T).$ Then we have the following Proposition 4.1.  \\\\
\textbf{Proposition 4.1} For $n\geq1$ and $t\in[0,T],$ we have $$Y_t\geq y_t^{n+1}\geq y_t^n,\ \ \ P-a.s.$$
\emph{Proof.} With the help of Lemma 3.7, Remark 4.1, Lemma 3.3 and Theorem 4.2, we can obtain this proposition from the argument of Coquet et al. [5, Lemma 6.2] or Hu et al. [9, Lemma 5.3], immediately. \ \ $\Box$\\

Set $$A_t^n:=\int_0^tn(Y_s-y_s^n)ds,\eqno(32)$$
Clearly, $A_t^n$ belongs to $S^\infty_{{\cal{F}}}(0,T)$ and is increasing with $A_0^n=0.$ By (31) and (32), we get that $\forall t\in [0,T],$
$${\cal{E}}[Y_T+zB_T+A_T^n-A_t^n|{\cal{F}}_t]=y_t^n+zB_t.\eqno(33)$$
By (H2), we have $\forall t\in [0,T],$
$${\cal{E}}[Y_T+zB_T+A_T^n|{\cal{F}}_t]=y_t^n+zB_t+A_t^n.$$
Thus $y_t^n+zB_t+A_t^n$ is an ${\cal{E}}$-martingale, by Lemma 3.4, there exists a pair $(g^n_s, Z^n_s)$  in $L^2_{{\cal{F}}}(0,T)\times L^2_{{\cal{F}}}(0,T;{\textbf{R}}^d)$ such that
$$g_s^n-g_s^m\leq\phi(|Z^n_s-Z^m_s|),\ \ \forall m\geq1,\eqno(34)$$
and$$y_t^n+zB_t+A_t^n=Y_T+zB_T+A_T^n+\int_t^Tg^n_sds-\int_t^TZ^n_sdB_s.$$
Then
$$y_t^n=Y_T+A_T^n-A_t^n+\int_t^Tg^n_sds-\int_t^T(Z^n_s-z)dB_s.\eqno(35)$$
Now, we can get\\\\
\textbf{Proposition 4.2} There exists a constant $C$ independent on $n$, such that
$$\textrm{(i)}\ \ E\int_0^T|Z^n_s-z|^2ds\leq C\ \ \textrm{ and }\ \ \textrm{(ii)} \ \ E|A^n_T|^2\leq C.$$
\emph{Proof.} In this proof, $C$  is assumed as a constant independent on $n$, its value may change line by line. By Proposition 4.1, we get that $y_t^1\leq y_t^n\leq y_t^{n+1}\leq Y_t.$ Thus, we have $$\|y^n_t\|_{L_{\cal{F}}^\infty}\leq C,\ \ n\geq1.\eqno(36)$$
By (35), (36), (9) and the fact that $\phi$ is increasing, subadditive and has a linear growth, we have
\begin{eqnarray*}
E|A_T^n|^2&\leq& 3E|y^n_0-y^n_T|^2+3TE\int_0^T|g^n_s|^2ds+3E\int_0^T|Z^n_s-z|^2ds\\
&\leq&C+ 3TE\int_0^T|\phi(|Z_s^n-z|)+\phi(|z|)|^2ds+3E\int_0^T|z^n_s-z|^2ds\\
&\leq&C+ 3TE\int_0^T(4\nu^2|Z_s^n-z|^2+4\nu^2+2(\phi(|z|))^2)ds+3E\int_0^T|Z^n_s-z|^2ds\\
&\leq&C+ 3(4\nu^2T+1)E\int_0^T|Z_s^n-z|^2ds\\
\end{eqnarray*}
Applying It\^{o} formula to $|y^n_t|^2,$ and by (36), (9), the fact that $\phi$ is increasing, subadditive and has a linear growth, and the inequality $2ab\leq \beta a^2+\frac{b^2}{\beta},\ \beta>0,$ we have
\begin{eqnarray*}|y_0^n|^2+E\int_0^T|Z^n_s-z|^2ds&=&E|Y_T|^2+2E\int_0^Ty_s^ng^n_sds
+2E\int_0^Ty_s^ndA_s^n\\
&\leq& C+2E\int_0^T|y_s^n||\phi(|Z_s^n-z|)+\phi(|z|)|ds
+2E\int_0^T|y_s^n|dA_s^n\\
&\leq& C+2E\int_0^T|y_s^n|(\nu|Z_s^n-z|+\nu+\phi(|z|))ds
+C[E|A_T^n|^2]^\frac{1}{2}\\
&\leq& C+\frac{1}{4}E\int_0^T|Z_s^n-z|^2ds+\frac{1}{6(4\nu^2T+1)}E|A_T^n|^2\\
\end{eqnarray*}
By above two inequalities, we can complete the proof. \ \ $\Box$\\\\

By (32), Proposition 4.1 and (ii) in Proposition 4.2, we get that as $n\rightarrow\infty,$
$$y^n_t\nearrow Y_t,\ \ dP\times dt-a.e.\eqno(37)$$
Then by (36) and Lebesgue dominated convergence theorem, we have
$$y^n_t\rightarrow Y_t,\ \ dt-a.e.\eqno(38)$$
in $L^2({\cal{F}}_T)$ sense.
By (9), (i) in Proposition 4.2 and linear growth of $\phi$, there exists a constant $C$ independent $n$ such that
$$E\int_0^T|g^n_s|^2ds\leq C.\eqno(39)$$
With the help of (36)-(39) and Proposition 4.2, we can apply the monotonic limit theorem (see Peng [18, Theorem 2.1] or Peng [19, Theorem 7.2]) to (35), then we have
$$Y_t=Y_T+A_T-A_t+\int_t^Tg_sds-\int_t^T(Z_s-z)dB_s.\eqno(40)$$
where $(Z_s-z)\in L_{\cal{F}}^2(0,T,\textbf{R}^d)$, $g_s\in L_{\cal{F}}^2(0,T)$ are the weak limits of $Z_s^n-z$ and $g_s^n$ in $L_{\cal{F}}^2(0,T,\textbf{R}^d)$ and $L_{\cal{F}}^2(0,T),$ respectively. For $t\in[0,T],\ A_t$ is the weak limit of $A_t^n$ in $L^2({\cal{F}}_T)$ and $A_t\in{\mathcal{D}}^2_{\cal{F}}(0,T)$ is increasing with $A_0=0$. Since $Y_t$ is continuous, then by (40), $A_t$ is a continuous process and by the monotonic limit theorem in Peng [18, 19] again, we further have
$$Z_t^n-z\rightarrow Z_t-z,\eqno(41)$$
in $L_{\cal{F}}^2$ sense, as $n\rightarrow \infty$. Then by (34),  (41) and the fact that $\phi$ is continuous with $\phi(0)=0,$ we can deduce that there exists a subsequence we still denote by ${n},$ such that the limit of $g_t^n$ exists, $dP\times dt-a.e.$ Thus by (39), Lebesgue dominated convergence theorem and the fact that $g_s\in L_{\cal{F}}^2(0,T)$ is the weak limit of $g_s^n$ in $L_{\cal{F}}^2(0,T),$ we can get
$$g_t^n\rightarrow g_t,\eqno(42)$$
in $L_{\cal{F}}^2$ sense, as $n\rightarrow \infty$. Thanks to  (38), (41) and (42), by (35) and (40), we can deduce that
$$A_T^n-A_t^n\rightarrow A_T-A_t,\ \  dt-a.e.\eqno(43)$$
in $L^2({\cal{F}}_T)$ sense, as $n\rightarrow \infty$. Then by (33), (38), (43) and (v) in Lemma 3.1, we can deduce that
$${\cal{E}}[Y_T+zB_T+A_T-A_t|{\cal{F}}_t]=Y_t+zB_t,\ \ dP\times dt-a.e.$$
By the continuity $Y_t$ and (H2), we have $\forall t\in[0,T],$
$${\cal{E}}[Y_T+zB_T+A_T|{\cal{F}}_t]=Y_t+zB_t+A_t,\ \ P-a.s.$$
The proof is complete. \ \ $\Box$
\section{Representation for ${\cal{F}}$-expectation by $g$-expectation}
The following representation theorem is the main result of this paper.\\\\
\textbf{Theorem 5.1} Let ${\cal{F}}$-expectation ${\cal{E}}$ satisfy (H1),  Then there exists a function $g(t,z): \Omega \times [0,T]\times {\mathbf{
R}}^d\longmapsto \mathbf{R}$ satisfying (A1), (A2) and (A3), such that, for each $X\in L^2({\cal{F}}_T)$ and $t\in[0,T],$ we have
\begin{center}
${\cal{E}}[X|{\cal{F}}_t]={\cal{E}}^g[X|{\cal{F}}_t],\ \ P-a.s.$
\end{center}
\emph{Proof.}
For $z\in \textbf{R}^d$, we consider the following SDE:
$$dY_t^z=-\phi(|z|)dt+zdB_t,\ \ Y_0^z=0.$$
Then
$$Y_t^z=Y_T^z+\int_t^T\phi(|z|)ds-\int_t^TzdB_s,\ \ Y_0^z=0.\eqno(44)$$
Clearly, $-\phi(|z|)t+zB_t=Y_t^z$ is a $\phi$-martingale and $-\phi(|z|)t\in S^\infty_{{\cal{F}}}(0,T).$ Then by (iii) in Lemma 3.1, we can check that $-\phi(|z|)t+zB_t$ is an ${\cal{E}}$-supermartingale. From the Theorem 4.3, there exists a process $A_t^z\in S^2_{{\cal{F}}}(0,T)$, which is increasing with $A_0^z=0$ such that
$${\cal{E}}[-\phi(|z|)T+zB_T+A_T^z|{\cal{F}}_t]=-\phi(|z|)t+zB_t+A_t^z,\ \ \ \forall t\in [0,T].\eqno(45)$$
Then by Lemma 3.4,  there exists a pair $(g(s,z), Z^z_s)$ such that
$$-\phi(|z|)t+zB_t+A_t^z=-\phi(|z|)T+zB_T+A_T^z+\int_t^Tg(s,z)ds-\int_t^TZ_s^zdB_s.\eqno(46)$$
and
$$|g(s,z)|\leq\phi(|Z_s^z|) \ \ \textrm{and}\ \ |g(s,z)-g(s,\bar{z})|\leq \phi(|Z_s^z-Z_s^{\bar{z}}|),\ \ \textrm{for}\ \ \bar{z}\in \textbf{R}^d.\eqno(47)$$
Comparing the bounded variation parts and martingale parts and of (44) and (46), we get
\begin{eqnarray*}
\phi(|z|)t&\equiv&A_t^z+\int_0^tg(s,z)ds,\\
 Z_s^z&\equiv& z.
\end{eqnarray*}
Then combining above equalities and (46), (47), we have
$$-\phi(|z|)t+zB_t+A_t^z=-\phi(|z|)T+zB_T+A_T^z+\int_t^Tg(s,z)ds-\int_t^TzdB_s.\eqno(48)$$
and
$$|g(s,z)|\leq\phi(|z|)\ \ \textrm{ and}\ \ |g(s,z)-g(s,\bar{z})|\leq  \phi(|z-\bar{z}|).\eqno(49)$$
Thus $g(t,z)$ satisfies (A1), (A2) and (A3).
By (48), (H2) and (45), we can get  for $0\leq r \leq t\leq T$
\begin{eqnarray*}
{\cal{E}}\left[-\int_r^tg(s,z)ds+\int_r^tzdB_s|{\cal{F}}_r\right]
&=&{\cal{E}}\left[-\phi(|z|)t+zB_t+A_t^z-(-\phi(|z|)r+zB_r+A_r^z)|{\cal{F}}_r\right]\\
&=&{\cal{E}}\left[-\phi(|z|)t+zB_t+A_t^z|{\cal{F}}_r\right]-(-\phi(|z|)r+zB_r+A_r^z)\\
&=&0.
\end{eqnarray*}
Thanks to the above equality, (49), the fact that $\phi$ is continuous with $\phi(0)=0$ and has a linear growth and (v) in Lemma 3.1, using the same argument as (7.4) in Coquet et al. [5], we can get for  $0\leq r \leq t\leq T$ and $\forall \eta_s\in L^2_{{\cal{F}}}(0,T,\textbf{R}^d),$
$${\cal{E}}\left[-\int_r^tg(s,\eta_s)ds+\int_r^t\eta_sdB_s|{\cal{F}}_r\right]=0.\eqno(50)$$
For $X \in L^2({\cal{F}}_T),$ we consider the BSDE
$$Y_t=X+\int_t^Tg(s,Z_s)ds-\int_t^TZ_sdB_s.\eqno(51)$$
By (49), BSDE (51) has a unique solution $(Y_t,Z_t)$ in ${\mathcal{S}}^2_{\cal{F}}(0,T)\times L^2_{\cal{F}}(0,T;{\mathbf{R}}^d).$ By (49) and Definition 2.1, we have
$${\cal{E}}^g\left[X|{\cal{F}}_t\right]=Y_t.$$
By (H2), (51) and (50), we get
$${\cal{E}}\left[X|{\cal{F}}_t\right]-Y_t={\cal{E}}\left[X-Y_t|{\cal{F}}_t\right]\\
={\cal{E}}\left[-\int_t^Tg(s,Z_s)ds+\int_t^TZ_sdB_s|{\cal{F}}_t\right]=0.$$
From above two equalities, we have
$${\cal{E}}[X|{\cal{F}}_t]={\cal{E}}^g[X|{\cal{F}}_t],\ \ P-a.s.$$

Now, we prove the uniqueness of $g$. Suppose there exists  another function $\bar{g}(t,z): \Omega \times [0,T]\times {\mathbf{
R}}^d\longmapsto \mathbf{R}$ satisfying (A1), (A2) and (A3), such that for each $X\in L^2({\cal{F}}_T)$ and $t\in[0,T],$ we have
$${\cal{E}}^g[X|{\cal{F}}_t]={\cal{E}}^{\bar{g}}[X|{\cal{F}}_t],\ \ P-a.s.$$
For each $z\in \textbf{R}^d$, $t\in[0,T]$ and $\varepsilon\in[0,T-t]$, Let $y_s^{t+\varepsilon}$ and $\bar{y}_s^{t+\varepsilon}$ are solutions of BSDEs with parameters $(g,z(B_{t+\varepsilon}-B_t),t+\varepsilon)$ and $(\bar{g},z(B_{t+\varepsilon}-B_t),t+\varepsilon),$ respectively. By (A3), we can check that
$$y_s^{t+\varepsilon}={\cal{E}}^g[z(B_{t+\varepsilon}-B_t)|{\cal{F}}_s]\ \ \textrm{and}\ \ \bar{y}_s^{t+\varepsilon}={\cal{E}}^{\bar{g}}[z(B_{t+\varepsilon}-B_t)|{\cal{F}}_s],\ \ s\in[0,t+\varepsilon].$$
Thus we have
$$y_s^{t+\varepsilon}=\bar{y}_s^{t+\varepsilon},\ \  P-a.s.,\ \ s\in[0,t+\varepsilon].$$
Then by the representation theorem for generator of BSDEs with continuous and linear growth generators (see Jia [11, Theorem 3.4] or Fan and Jiang [7, Theorem 2]) and some simple arguments, we can get that $\forall z\in\textbf{R}^d,$
$$g(t,z)=\bar{g}(t,z),\ \ dP\times dt-a.e.$$
By (A1), we have $dP\times dt-a.e.,$
$$g(t,z)=\bar{g}(t,z),\ \ z\in\textbf{R}^d.$$
 The proof is complete. \ \ $\Box$\\

 Theorem 4.1 and Theorem 4.2 are existence and uniqueness theorem and comparison theorem of ${\cal{E}}(f,T,X,z),$  respectively, with $X\in L^\infty({\cal{F}}_T)$ and $f(\cdot,0)\in L^\infty_{{\cal{F}}}(0,T)$. From Theorem 5.1, we can get the following general result.\\\\
 \
\textbf{Corollary 5.1} Let ${\cal{F}}$-expectation ${\cal{E}}$ satisfy (H1), $z\in\textbf{R}^d,$ $X\in L^2({\cal{F}}_T)$ and $f(\cdot,0)\in L^2_{{\cal{F}}}(0,T).$ Then ${\cal{E}}(f,T,X,z)$ has a unique solution $y_t\in {\cal{S}}_{\cal{F}}^2(0,T).$ Moreover if $\bar{y}_t$ is the solution of the following ${\cal{E}}(f+\eta_t,T,\bar{X},z)$:
$$\bar{y}_t+zB_t={\cal{E}}\left[\bar{X}+zB_T+\int_t^T(f(s,\bar{y}_s)+\eta_s)ds|{\cal{F}}_t\right],\ \  t\in[0,T],$$
where $\bar{X}$ in $L^2({\cal{F}}_T)$ and $\eta_t\in L^2_{{\cal{F}}}(0,T)$ satisfy
$$\bar{X}\geq X,\ \ \ \eta_t\geq0,\ \ dP\times dt-a.e.,$$
then $\forall t\in[0,T],$ we have
$$\bar{y}_t\geq y_t,\ \ \ P-a.s.$$
\emph{Proof.} By the proof of Theorem 5.1,  there exist a function $g: \Omega \times [0,T]\times {\mathbf{
R}}^d\longmapsto \mathbf{R}$ satisfying (A1), (A2) and (A3) such that $g$ satisfies (50) and for $\xi$ in $L^2({\cal{F}}_T),$ $${\cal{E}}[\xi|{\cal{F}}_t]={\cal{E}}^g[\xi|{\cal{F}}_t],\ \ P-a.s.\eqno(52)$$
Let $z\in\textbf{R}^d,$ $X\in L^2({\cal{F}}_T)$ and $f(\cdot,0)\in L^2_{{\cal{F}}}(0,T).$ Set $\hat{g}(t,\hat{y},\hat{z}):=f(t,\hat{y}-zB_t)+g(t,\hat{z}).$ Clearly, $\hat{g}(t,\hat{y},\hat{z})$ satisfies (A1) and (A2), Thus, the BSDE
$$Y_t=X+zB_T+\int_t^Tf(s,Y_s-zB_s)ds+\int_t^Tg(s,Z_s)ds-\int_t^TZ_sdB_s.\eqno(53)$$
has a unique solution $(Y_t,Z_t)\in{\mathcal{S}}^2_{\cal{F}}(0,T)\times L^2_{\cal{F}}(0,T;{\mathbf{R}}^d)$.
By (53), (50) and (H2), we can get
$$Y_t={\cal{E}}\left[X+zB_T+\int_t^Tf(s,Y_s-zB_s)ds|{\cal{F}}_t\right].$$
Then by setting $y_t:={Y}_t-zB_t$, we get ${\cal{E}}(f,T,X,z)$ has a solution $y_t\in S^2_{{\cal{F}}}(0,T).$ Moreover, by (52) and the uniqueness of solution of BSDE (53), we also can deduce $y_t$ is a unique solution. In fact, if ${\cal{E}}(f,T,X,z)$ has another solution $\hat{y}_t\in S^2_{{\cal{F}}}(0,T),$ then there exists a process $\hat{Z}_t\in L^2_{\cal{F}}(0,T;{\mathbf{R}}^d)$ such that
\begin{eqnarray*}
\hat{y}_t+zB_t&=&{\cal{E}}\left[X+zB_T+\int_t^Tf(s,\hat{y}_s)ds|{\cal{F}}_t\right]\\
&=&{\cal{E}}\left[X+zB_T+\int_0^Tf(s,\hat{y}_s)ds|{\cal{F}}_t\right]-\int_0^tf(s,\hat{y}_s)ds\\
&=&{\cal{E}}^g\left[X+zB_T+\int_0^Tf(s,\hat{y}_s)ds|{\cal{F}}_t\right]-\int_0^tf(s,\hat{y}_s)ds\\
&=&X+zB_T+\int_t^Tf(s,\hat{y}_s)ds+\int_t^Tg(s,\hat{Z}_s)ds-\int_t^T\hat{Z}_sdB_s.
\end{eqnarray*}
Thus, by the uniqueness of solution of BSDE (53), we can get $dP\times dt-a.e.,\ y_t=\hat{y}_t.$

Set $\bar{g}(t,\bar{y},\bar{z}):=f(t,\bar{y}-zB_t)+\eta_t+g(t,\bar{z}).$ Clearly, $\bar{g}(t,\bar{y},\bar{z})$ also satisfies (A1) and (A2).  Thus, the BSDE
$$\bar{Y}_t=\bar{X}+zB_T+\int_t^T(f(s,\bar{Y}_s-zB_s)+\eta_s)ds+\int_t^Tg(s,\bar{Z}_s)ds-\int_t^T\bar{Z}_sdB_s.\eqno(54)$$
has a unique solution $(\bar{Y}_t,\bar{Z}_t)\in{\mathcal{S}}^2_{\cal{F}}(0,T)\times L^2_{\cal{F}}(0,T;{\mathbf{R}}^d)$. By the argument as above, we get $\bar{y}_t:={\bar{Y}}_t-zB_t$ is the unique solution of ${\cal{E}}(f+\eta_t,T,\bar{X},z).$ Then by (53), (54) and comparison theorem for BSDEs under (A1) and (A2) (see Jia [11, Theorem 3.1]), we can get $\forall t\in[0,T],$ $\bar{Y}_t\geq Y_t, \ P-a.s.$ Thus $\forall t\in[0,T],$ $\bar{y}_t\geq y_t, \ P-a.s.$ The proof is complete. \ \ $\Box$

In fact, if ${\cal{F}}$-expectation ${\cal{E}}$ satisfy (H1), by Theorem 5.1, Lemma 2.1 and similar argument as Corollary 5.1, we also can get a Doob-Meyer decomposition for ${\cal{E}}$-supermartingale $y_t+zB_t$ with $y_t\in S^2_{{\cal{F}}}(0,T).$ We leave it to the interested readers.

\end{document}